\input amstex
\documentstyle{amsppt}
\magnification=\magstep1
 \hsize 13cm \vsize 18.35cm \pageno=1
\loadbold \loadmsam
    \loadmsbm
    \UseAMSsymbols
\topmatter
\NoRunningHeads
\title $q$-Bernoulli numbers and polynomials
associated with Gaussian binomial coefficient
\endtitle
\author
  Taekyun Kim
\endauthor

\abstract The first purpose of this paper is to present a systemic
 study of some families of multiple $q$-Bernoulli numbers and
 polynomials by using multivariate $q$-Volkenborn integral (= $p$-adic $q$-integral) on $\Bbb
 Z_p $. From the studies of these $q$-Bernoulli numbers and
 polynomials of higher order we derive some interesting $q$-analogs
 of Stirling number identities.

\endabstract
\thanks  2000 AMS Subject Classification: 11B68, 11S80
\newline  This paper is supported by  Jangjeon Research Institute for Mathematical
Science(JRIMS-10R-2001)
\endthanks
\endtopmatter

\document

{\bf\centerline {\S 1. Introduction}}

 \vskip 20pt

Let $q$ be regarded as either a complex number $q\in\Bbb C$ or a
$p$-adic number $q\in\Bbb C_p$. If $q\in\Bbb C$, then we always
assume $|q|<1$. If $q\in\Bbb C_p$, we normally assume
$|1-q|_p<p^{-\frac{1}{p-1}},$ which implies that $q^x =\exp(x\log q
) $ for $|x|_p\leq 1$. Here, $| \cdot |_p$ is the $p$-adic absolute
value in $\Bbb C_p$ with $|p|_p=\frac{1}{p}.$ The $q$-basic natural
number are defined by $[n]_q=\frac{1-q^n}{1-q}=1+q+\cdots+q^{n-1}, $
( $n\in \Bbb N $), and $q$-factorial are also defined  as $[n]_q
!=[n]_q \cdot [n-1]_q \cdots [2]_q \cdot [1]_q.$ In this paper we
use the notation of Gaussian binomial coefficient as follows:
$${\binom {n}{k}}_q =\frac{[n]_q!}{[n-k]_q![k]_q!}
=\frac{[n]_q\cdot [n-1]_q \cdots [n-k+1]_q}{[k]_q!}. \tag1$$ Note
that $\lim_{q\rightarrow 1}{\binom{n}{k}}_q=\binom{n}{k}=\frac{n
\cdot(n-1)\cdots(n-k+1)}{n!}.$ The Gaussian coefficient satisfies
the following recursion formula:
$$\binom{n+1}{k}_q =\binom{n}{k-1}_q+q^k\binom{n}{k}_q=q^{n-k}\binom{n}{k-1}_q+\binom{n}{k}_q, \text{ cf. [1-23].} \tag2$$
From thus recursion formula we derive
$$\binom{n}{k}_q=\sum_{d_0+\cdots+d_k=n-k, d_i\in\Bbb N} q^{d_1+2d_2+\cdots+kd_k},\text{ see [15, 20, 21].}\tag3$$
Let $p$ be a fixed prime. Throughout this paper $\Bbb Z_p ,$ $\Bbb
Q_p ,$ $\Bbb C,$ and $\Bbb C_p$ will, respectively, denote the ring
of $p$-adic rational integers, the field of $p$-adic rational
numbers, the complex number field, and the completion of algebraic
closure of $\Bbb Q_p .$ For $d$ a fixed positive integer $(p,d)=1$,
let
$$\split
& X=X_d = \lim_{\overleftarrow{N} } \Bbb Z/ dp^N \Bbb Z , \text{ and
 }X_1 = \Bbb Z_p , \cr  & X^\ast = \underset {{0<a<d p}\atop
{(a,p)=1}}\to {\cup} (a+ dp \Bbb Z_p ), \cr & a+d p^N \Bbb Z_p =\{
x\in X | x \equiv a \pmod{dp^N}\},\endsplit$$ where $a\in \Bbb Z$
lies in $0\leq a < d p^N$, cf. [8-18]. For $x\in\Bbb C_p ,$ we use
the notation $[x]_q=\frac{1-q^x}{1-q}, $ cf. [1-6].

We say that $f$ is a uniformly differentiable function at a point $a
\in\Bbb Z_p $ and denote this property by $f\in UD(\Bbb Z_p )$, if
the difference quotients $F_f (x,y) = \dfrac{f(x) -f(y)}{x-y} $ have
a limit $l=f^\prime (a)$ as $(x,y) \to (a,a)$. For $f\in UD(\Bbb Z_p
)$, let us start with the expression

$$\eqalignno{ & \dfrac{1}{[p^N ]_q} \sum_{0\leq j < p^N} q^j f(j) =\sum_{0\leq j < p^N} f(j)
\mu_q (j +p^N \Bbb Z_p ), }
$$
representing a $q$-analogue of Riemann sums for $f$, cf. [8, 21-23].
The integral of $f$ on $\Bbb Z_p$ will be defined as limit ($n \to
\infty$) of those sums, when it exists. The $q$-Volkenborn integral
(=$p$-adic $q$-integral) of the function $f\in UD(\Bbb Z_p )$ is
defined by
$$ I_q (f) =\int_{X }f(x) d\mu_q (x) = \lim_{N\to \infty}
\dfrac{1}{[dp^N ]_q} \sum_{0\leq x < dp^N} f(x) q^x ,\quad \text{see
[8] .} \tag4 $$ The Carlitz's $q$-Bernoulli numbers $\beta_{k,q}$
can be determined inductively by
$$\eqalignno{ &
\beta_{0,q} =1,\quad  q(q\beta +1)^k -\beta_{k,q} = \cases
1 & \text{if\   $k=1$}\\
0 & \text{if \  $k>1$,}
\endcases }
$$
with the usual convention of replacing $\beta^i$ by
$\beta_{i,q}$(see [2, 3, 24, 25]).

In [8], it was shown that the Carlitz's $q$-Bernoulli numbers can be
represented by $p$-adic $q$-integral on $\Bbb Z_p$ as follows:
$$\int_{\Bbb
Z_p}[x]_q^md\mu_q(x)=\int_{X}[x]_q^m d\mu_q(x)=\beta_{m,q}, \text{
$m\in\Bbb Z_{+}$.}\tag5$$ The $k$-th order factorial of the
$q$-number $[x]_q$, which is defined by
$$[x]_{k,q}=[x]_q\cdot[x-1]_q\cdots[x-k+1]_q=\frac{(1-q^x)(1-q^{x-1})\cdots(1-q^{x-k+1})}{(1-q)^k},
$$ is called $q$-factorial of $x$ of order $k$, cf.[15, 19, 20].
From this we note that $\binom{x}{k}_q=\frac{[x]_{k,q}}{[k]_q!},$
cf.[4-16]. The theory of $q$-number and the factorial of $q$-number
are applicable in the many areas related to mathematics ,
mathematical physics and probability. For example, we consider a
sequence of Bernoulli trials and assume that the conditionally
probability of success at the $n$-th trial, given that $k$ successes
occur before that trial varies geometrically with $n$ and $k$.
specifically, suppose that the probability of success at the
$n+1$-th trial, given that $k$ success occur up the $n$-th trial, is
given by
$$\lambda_{n,k}=q^{an+bk+c} (\in \Bbb R ), \text{ $k=0,1,2, \cdots, n,$ $n=0,1,2,\cdots, $} $$
with $a$, $ b$ and $c$ such that $0\leq \lambda_{n,k} \leq 1$. The
particular case $b=0$ corresponds to the assumption that the
probability of success at any trial depends only on the number of
previous trials, while the other particular case $a=0$ corresponds
to the assumption that the probability of success at only trial
depends only on the number of previous success. The purpose of this
paper is to present a systemic study of some families of multiple
$q$-Bernoulli numbers and polynomials by using multivariate
$q$-Volkenborn integral(=$p$-adic $q$-integral) on $\Bbb Z_p$. From
the studies of these $q$-Bernoulli numbers and polynomials we derive
some interesting $q$-analogs of Stirling number identities. That is,
the $q$-analogs of many classical Stirling number identities are
formulated and their interesting features are revealed in this
paper.

\vskip 20pt

{\bf\centerline {\S 2. $q$-Bernoulli numbers associated with
$q$-Stirling number identities}} \vskip 10pt

In this section we assume that $q\in\Bbb C_p$ with
$|1-q|_p<p^{-\frac{1}{p-1}}$. From the definition of $[x]_q$ we can
easily derive the following equation.
$$q^n [x-n]_q=\frac{q^n-1+1-q^x}{1-q}=[x]_q-[n]_q, \text{ and }
[-x]_q=\frac{1}{q^x}\frac{q^x-1}{1-q}=-\frac{1}{q^x}[x]_q. \tag6$$
Let $(Eh)(x)=h(x+1)$ be the shift operator.  Then we  consider the
$q$-difference operator as follows:
$$\Delta_q^n = \prod_{i=1}^{n}\left(E-q^{i-1}I \right), \text{ where $(Ih)(x)=h(x)$.}\tag7$$
From (6) and (7), we note that
$$ f(x)=\sum_{n\geq 0}\binom{x}{n}_q \Delta_q^nf(0), \tag8$$
where
$$\Delta_q^n
f(0)=\sum_{k=0}^{n}\binom{n}{k}_q(-1)^kq^{\binom{k}{2}}f(n-k).
\tag9$$ The $q$-Stirling number of the second kind is defined by
Carlitz as follows:
$$s_2(n,k,q)=\frac{q^{-\binom{k}{2}}}{[k]_q!}\sum_{j=0}^k
(-1)^jq^{\binom{j}{2}}\binom{k}{j}_q[k-j]_q^n , \text{ see
[3].}\tag10$$ By (9) and (10) we easily see that
$$s_2(n,k,q)=\frac{q^{-\binom{k}{2}}}{[k]_q!}\Delta_q^k0^n. \tag
11$$ From (11) we can also derive the following equation.
$$[x]^n=\sum_{k=0}^n\binom{x}{k}_q[k]_q!s_2(k,n-k,q)q^{\binom{k}{2}}=\sum_{k=0}^n[x]_{k,q}
\frac{q^{\binom{k}{2}_q-\binom{n-k}{2}_q}}{[n-k]_q!}\Delta_q^{n-k}0^k
 . \tag12$$ By (2), we easily see
that
$$\int_{\Bbb Z_p}\binom{x}{n}_q
d\mu_q(x)=\frac{(-1)^n}{[n+1]_q}q^{(n+1)-\binom{n+1}{2}}. \tag13$$
From (5), (12) and (13) we can derive the following theorem.
\proclaim{ Theorem 1} For $m\in\Bbb Z_{+}$, we have
$$ \beta_{m,q}=q\sum_{k=0}^m
\frac{[k]_q!}{[k+1]_q}(-1)^k s_{2}(k,n-k,q),\tag14
$$ where $\beta_{m,q}$ are $m$-th Carlitz $q$-Bernoulli numbers.
 \endproclaim
The $q$-Stirling numbers of the first kind is defined as
$$(1-q)^n[x]_{n,q}=\prod_{i=1}^n(1-q^{x-n+1}q^{i-1})
=\sum_{l=0}^n \binom{n}{l}_q
q^{\binom{l}{2}}(-1)^lq^{l(x-n+1)}.\tag15$$ It is easy to see that
$$q^{lx}=\left([x]_q(q-1)+1\right)^l=\sum_{m=0}^l\binom{l}{m}(q-1)^m[x]_q^m.
\tag16$$ From (16) we note that
$$\aligned
&\frac{1}{(1-q)^n}\sum_{l=0}^n\binom{n}{l}_q
q^{\binom{l}{2}}(-1)^lq^{l(x-n+1)}\\
&=\frac{1}{(1-q)^n}\sum_{l=0}^n\binom{n}{l}_q q^{\binom{l}{2}+l-ln}(-1)^l\sum_{m=0}^l\binom{l}{m}(q-1)^m[x]_q^m\\
&=\frac{1}{(1-q)^n}\sum_{m=0}^n(q-1)^m\left(\sum_{l=m}^n\binom{n}{l}_qq^{\binom{l}{2}-ln+l}
\binom{l}{m}(-1)^l\right)[x]_q^m .
\endaligned \tag17$$
By (12), (15) and (17) we obtain the following theorem.

\proclaim{ Theorem 2} For $n\in\Bbb Z_{+}$ we have
$$\beta_{n,q}=\sum_{l=0}^n s_{2}(l,
n-l,q)q^{\binom{l}{2}}\sum_{m=0}^l\frac{1}{(1-q)^{l-m}}\left(\sum_{i=m}^l\binom{l}{i}_q\binom{i}{m}
q^{\binom{i}{2}-il+i}(-1)^i\right)\beta_{m,q}.
$$
\endproclaim
In [3] Carlitz has given the following relation.
$$s_2(n,k,q)=(q-1)^{-k}\sum_{j=0}^k(-1)^{k-j}\binom{k+n}{k-j}\binom{j+n}{j}_q,
\tag18 $$ and
$$\binom{n}{k}_q=\sum_{j=0}^n\binom{n}{j}(q-1)^{j-k}s_2(k,j-k,q).$$
By simple calculation we easily see that
$$q^{nt}=\sum_{k=0}^n(q-1)^kq^{\binom{k}{2}}\binom{n}{k}_q[t]_{k,q}
=\sum_{m=0}^n\left(\sum_{k=m}^n(q-1)^k\binom{n}{k}_qs_1(k,m,q)\right)[t]_q^m.
\tag19$$ By using $p$-adic $q$-integral on $\Bbb Z_p$ we have
$$\int_{\Bbb
Z_p}q^{nt}d\mu_q(t)=\sum_{m=0}^n\binom{n}{m}(q-1)^m\beta_{m,q}.
\tag20$$ From (19) and (20) we derive
$$\binom{n}{m}=\sum_{k=m}^n(q-1)^{-m+k}\binom{n}{k}_qs_1(k,m,q).
\tag21$$ From the definition of the first kind Stirling number  we
note that
$$q^{\binom{n}{2}}\binom{x}{n}_q[n]_q!=[x]_{n,q}q^{\binom{n}{2}}
=\sum_{k=0}^ns_1(n,k,q)[x]_q^k.\tag22$$ By (13) and (22) we have
$$\frac{1}{[n+1]_q}=\frac{q^{-1}}{[n]_q!}\sum_{k=0}^n(-1)^{n-k}s_1(n,k,q)\beta_{k,q}.\tag23$$
By (14), (15), (17) and (23) we obtain the following theorem.
\proclaim{Theorem 3} For $n,j\in\Bbb Z_{+}$ we have
$$s_1(n,j,q)=\frac{q^{\binom{n}{2}}}{(q-1)^{n-j}}\sum_{k=j}^n(-1)^{n-k}q^{\binom{k+1}{2}-nk}\binom{n}{k}_q\binom{k}{j}.
\tag 24$$ Moreover,
$$\frac{1}{[n+1]_q}=\frac{q^{-1}}{[n]_q!}\sum_{k=0}^n(-1)^{n-k}s_1(n,k,q)\beta_{k,q}.$$
\endproclaim
\vskip 10pt

 {\bf\centerline {\S 3. Multivariate $p$-adic $q$-integral on $\Bbb Z_p$ associated with $q$-Stirling numbers  }}
  \vskip 10pt

In this section we also assume that $q\in\Bbb C_p$ with
$|1-q|_p<p^{-\frac{1}{p-1}}.$ For any positive integers $k,m ,$ we
consider the following multivariate $p$-adic $q$-integral on $\Bbb
Z_p$ related to $q$-Bernoulli polynomials of higher order as
follows:
$$\beta_{n,q}^{(k)}(x)=\frac{1}{(1-q)^n}\sum_{i=0}^n(-1)^i\binom{n}{i}q^{ix}\int_{\Bbb
Z_p}\cdots\int_{\Bbb Z_p}q^{\sum_{l=1}^k(k-l+i)x_l}d\mu_q(x_1)\cdots
d\mu_{q}(x_k).\tag25$$ In the special case $x=0,$
$\beta_{n,q}^{(k)}(0)=\beta_{n,q}^{(k)}$ will be called the
$q$-Bernoulli numbers of order $k$. From (25) we note that
$$\beta_{n,q}^{(k)}(x)=\frac{1}{(1-q)^n}\sum_{i=0}^n(-1)^i
\binom{n}{i}\frac{(i+k)\cdots(i+1)}{[i+k]_q\cdots[i+1]_q}q^{ix}.\tag26$$
Thus, we obtain the following theorem.

\proclaim{Theorem 4} For $m,k\in\Bbb Z_{+},$ we have
$$\beta_{n,q}^{(k)}(x)=\frac{1}{(1-q)^n}\sum_{i=0}^n\frac{(-1)^i\binom{n}{i}\binom{i+k}{k}}{\binom{i+k}{k}_q}
\frac{k!}{[k]_q!}q^{ix}.$$
\endproclaim
Now we also define $\beta_{n,q}^{(-k)}(x)$ as follows:
$$\beta_{n,q}^{(-k)}(x)=\frac{1}{(1-q)^n}\sum_{i=0}^n
\frac{(-1)^i\binom{n}{i}q^{ix}}{\int_{\Bbb Z_p}\cdots\int_{\Bbb
Z_p}q^{\sum_{l=1}^k(k-l+i)x_l}d\mu_q(x_1)\cdots d\mu_q(x_k)},
\tag27$$ where $n,k$ are positive integers. From (27) we note that
$$\beta_{n,q}^{(-k)}(x)=\frac{1}{(1-q)^n}\sum_{i=0}^n(-1)^i\binom{n}{i}\frac{\binom{i+k}{k}_q}{\binom{i+k}{k}}
\frac{[k]_q!}{k!}q^{ix}.\tag28$$ It is easy to see that

$$\frac{\binom{k}{j}}{\binom{j+n}{n}
n!}=\frac{(k+n)\cdots(k+1)k\cdots(k-j+1)}{(j+n)!(k+n)\cdots(k+1)}=\frac{\binom{k+n}{k-j}}{\binom{k+n}{n}n!}.
\tag29$$ By (27), (28) and (29) we obtain the following theorem.
\proclaim {Theorem 5} For $n,k \in \Bbb Z_{+},$ we have
$$\beta_{n,q}^{(-k)}(x)=\frac{1}{(1-q)^n}\sum_{i=0}^n(-1)^i\binom{i+k}{k}_q\frac{\binom{n+k}{n-i}}{\binom{n+k}{k}}
\frac{[k]_q!}{k!}q^{ix}.\tag 30$$
 \endproclaim
From (18) and (30) we derive
$$s_2(n,k,q)=\binom{k+n}{n}\frac{n!}{[n]_q!}\beta_{k,q}^{(-k)}(0).$$
That is,
$$\frac{1}{(1-q)^k}\sum_{i=0}^k\frac{(-1)^k\binom{k}{i}}{\int_{\Bbb
Z_p}\cdots\int_{\Bbb Z_p}q^{\sum_{l=1}^n(n-l+i)x_l}d\mu_q(x_1)\cdots
d\mu_q(x_n)} =\frac{[n]_q!}{\binom{k+n}{n}n!}s_2(n,k,q).$$ Thus, we
note that
$$\beta_{0,q}^{(-k)}(0)=\left(\int_{\Bbb Z_p}\cdots\int_{\Bbb
Z_p}q^{\sum_{i=1}^k(k-i)x_i}d\mu_q(x_1)\cdots d\mu_q(x_k)
\right)^{-1}=\frac{[k]_q!}{k!}. \tag31$$ By the same method we see
that $\beta_{1,q}^{(2)}(0)=\frac{-2(q+2)} {[2]_q[3]_q} , \cdots,$
$\beta_{0,q}^{(-k)}(0)=\frac{[k]_q!}{k!}.$ Thus, we have
$$s_2(k,0,q)=\frac{k!}{[k]_q!}\beta_{0,q}^{(-k)}(0)=\frac{k!}{[k]_q!}\frac{[k]_q!}{k!}=1.$$
From the definition of $\beta_{m,q}^{(k)}$ we can also derive the
following equality.
$$\aligned
&\sum_{i=0}^m\binom{m}{i}(q-1)^i\int_{\Bbb Z_p}\cdots\int_{\Bbb
Z_p}[x_1+\cdots+x_k]_q^iq^{\sum_{l=1}^k(k-l)x_l}d\mu_q(x_1)\cdots
d\mu_q(x_k)\\
&=\int_{\Bbb Z_p}\cdots\int_{\Bbb
Z_p}q^{(m+k-1)x_1+\cdots+(m+1)x_{k-1}+mx_k}d\mu_q(x_1)\cdots
d\mu_q(x_k)=\frac{\binom{m+k}{k}}{\binom{m+k}{k}_q}\frac{k!}{[k]_q!}.
\endaligned$$
Therefore we obtain the following:
$$\sum_{i=0}^m\binom{m}{i}(q-1)^i\beta_{i,q}^{(k)}=\frac{\binom{m+k}{k}}{\binom{m+k}{k}_q}
\frac{k!}{[k]_q!}.$$ Finally, we observe that
$$q^{\binom{n}{2}}[x]_{n,q}=[x]_q q[x-1]_q\cdots q^{n-1}[x-n+1]_q
=[x]_q\cdot([x]_q-1)\cdots([x]_q-[n-1]_q).$$ Thus, we have
$$q^{\binom{n}{2}}\binom{x}{n}_q=\frac{1}{[n]_q!}\prod_{k=0}^n \left([x]_q-[k]_q \right)
=\frac{1}{[n]_q!}\sum_{k=0}^ns_1(n,k,q)[x]_q^k.$$
\vskip 10pt

 {\bf\centerline {\S 4. Further Remarks and Observations  }}
  \vskip 10pt
In this section, let $p$ be a fixed odd prime number. For $n\in \Bbb
N$, $k\in\Bbb Z_{+} ,$ we consider the following $q$-Euler numbers
of higher order.
$$E_k^{(n)}(x,q)=\int_{\Bbb Z_p}\cdots \int_{\Bbb
Z_p}[\sum_{i=1}^nx_i +x
]_q^kq^{\sum_{j=1}^nx_j(n-j)}d\mu_{-q}(x_1)\cdots d\mu_{-q}(x_n),
\tag32$$ where $\mu_{-q}(x+p^N \Bbb
Z_p)=\frac{1+q}{1+q^{p^N}}(-q)^x=\frac{(-q)^x}{[p^N]_{-q}}$, see
[10]. From (32) we note that
$$E_n^{(k)}(x,q)=\frac{[2]_q^n}{(1-q)^k}
\sum_{l=0}^k\binom{k}{l}\frac{(-1)^lq^{lx}}{(1+q^{n+l})\cdots(1+q^{l+1})}.$$
The $q$-binomial formulae are known as
$$\prod_{i=1}^n\left(a+bq^{i-1}\right)=\sum_{k=0}^n\binom{n}{k}_qq^{\binom{k}{2}}q^{n-k}b^k
,\tag33$$ and
$$\prod_{i=1}^n\left(a-bq^{i-1}\right)^{-1}=\sum_{k=0}^{\infty}\binom{n+k-1}{k}_q
b^k. $$ By (32) and (33) we obtain the following:

\proclaim{ Proposition 6} For $n\in\Bbb N$, $k\in\Bbb Z_{+}$, we
have
$$E_k^{(n)}(x,q)=\frac{[2]_q^n}{(1-q)^k}\sum_{l=0}^k\binom{k}{l}(-1)^lq^{lx}
\sum_{i=0}^{\infty}\binom{n+i-1}{i}_q(-1)^i q^{(l+1)i}.$$
\endproclaim
 Seeking to define  a suitable polynomial analogue for negative
 value of $n$, we give the definition as follows:
 $$E_k^{(-n)}(x,q)=\frac{1}{(1-q)^k}\frac{1}{[2]_q^n}
 \sum_{l=0}^k\binom{k}{l}(-1)^lq^{lx}\left(\prod_{i=1}^n(1+q^{l+i})\right).\tag34$$
From (33) and (34) we can also derive the following Eq.(35).
$$E_k^{(-n)}(x,q)=\frac{1}{(1-q)^k}\frac{1}{[2]_q^n}\sum_{l=0}^k\binom{k}{l}(-1)^lq^{lx}
\sum_{i=0}^n\binom{n}{i}_qq^{\binom{i}{2}}q^{(l+1)i}.\tag35$$

\vskip 20pt

\quad Taekyun Kim

\quad EECS, Kyungpook National University, Taegu 702-701, S. Korea

\quad e-mail:\text{ tkim$\@$knu.ac.kr; tkim64$\@$hanmail.net}

\enddocument